\documentclass[10pt,reqno]{amsart}

\usepackage{amssymb}
\usepackage{amsmath}
\usepackage{amsfonts}
\usepackage{amsthm}
\usepackage{paralist}

\usepackage{mathrsfs}

\usepackage{graphicx,multirow,array}

\numberwithin{equation}{section}

\usepackage{subfigure}
\usepackage{epstopdf}

\usepackage{ifthen} 

\provideboolean{shownotes} 
\setboolean{shownotes}{true} 
\newcommand{\margnote}[1]{
\ifthenelse{\boolean{shownotes}}%
{\marginpar{\raggedright\tiny\texttt{#1}}}%
{}%
}

\newcommand{\hole}[1]{
\ifthenelse{\boolean{shownotes}}%
{\begin{center} \fbox{ \rule {.25cm}{0cm}
\rule[-.1cm]{0cm}{.4cm} \parbox{.85\textwidth}{\begin{center}
\texttt{#1}\end{center}} \rule {.25cm}{0cm}}\end{center}}
{}
}

\theoremstyle{plain}

\theoremstyle{definition}

\theoremstyle{remark}

\newcommand{\prt}{\partial_t}
\newcommand{\prtt}{\partial_{tt}}
\newcommand{\prx}{\partial_x}
\newcommand{\prxx}{\partial_{xx}}
\newcommand{\flux}{\mathcal{A}}
\newcommand{\spgw}{\mathcal{S}}
\newcommand{\diag}{\mathcal{D}}
\newcommand{\R}{\mathbb{R}}
\newcommand{\Z}{\mathbb{Z}}
\newcommand{\N}{\mathbb{N}}
\newcommand{\xx}{{\rm x}}
\newcommand{\uu}{{\rm u}}
\newcommand{\vv}{{\rm v}}

\newcommand{\xim}{{\rm x}_{i-\frac12}}
\newcommand{\xip}{{\rm x}_{i+\frac12}}
\newcommand{\dxi}{{\rm dx}_i}
\newcommand{\dxim}{{\rm dx}_{i-1}}
\newcommand{\dxip}{{\rm dx}_{i+1}}
\newcommand{\dtn}{{\rm dt}_n}
\newcommand{\sd}{{\rm s}}
\newcommand{\rd}{{\rm r}}

\newcommand{\blankbox}[2]{%
  \parbox{\columnwidth}{\centering
    \setlength{\fboxsep}{0pt}%
    \fbox{\raisebox{0pt}[#2]{\hspace{#1}}}%
  }%
}

\begin{document}

\title[Kinetic schemes for assessing stability of traveling fronts]{Kinetic schemes for assessing stability of traveling fronts for the Allen-Cahn equation with relaxation}

\author[C. Lattanzio]{Corrado Lattanzio}

\address{{\rm (C. Lattanzio)} Dipartimento di Ingegneria e Scienze dell'Informazione e Matematica\\Universit\`a degli Studi dell'Aquila\\via Vetoio (snc), Coppito I-67010, L'Aquila (Italy)}

\email{corrado@univaq.it}

\author[C. Mascia]{Corrado Mascia}

\address{{\rm (C. Mascia)} Dipartimento di Matematica `G. Castelnuovo'\\Universit\'a di Roma `La Sapienza'\\Pia\-zzale A. Moro 2, I-00185 Roma (Italy)}

\email{mascia@mat.uniroma1.it}

\author[R.G. Plaza]{Ram\'on G. Plaza} 

\address{{\rm (R. G. Plaza)} Instituto de Investigaciones en Matem\'aticas Aplicadas y en 
Sistemas\\Universidad Nacional Aut\'onoma de M\'exico\\Circuito Escolar s/n C.P. 04510 Cd. de M\'exico 
(Mexico)}

\email{plaza@mym.iimas.unam.mx}

\author[C. Simeoni]{Chiara Simeoni}

\address{{\rm (C. Simeoni)} Laboratoire J.A. Dieudonn\'e UMR CNRS 7351\\Universit\'e de Nice 
Sophia-Antipolis\\Parc Valrose 06108 Nice Cedex 02 (France)}

\email{simeoni@unice.fr}



\begin{abstract}
This paper deals with the numerical (finite volume) approximation of reaction-diffusion systems with relaxation, among which the hyperbolic extension of the Allen--Cahn equation represents a notable prototype. Appropriate discretizations are constructed starting from the kinetic interpretation of the model as a particular case of reactive jump process. Numerical experiments$^*$\protect{\thanks{$^*$the code for reproducing the numerical tests is available upon request to the authors.}} are provided for exemplifying the theoretical analysis (previously developed by the same authors) concerning the stability of traveling waves, and important evidence of the validity of those results beyond the formal hypotheses is numerically established.
\end{abstract}

\maketitle

\setcounter{tocdepth}{1}
%


\section{Physical motivations and problem statement}
\label{sec:physics}

The standard approach to heat conduction in a medium is based on the continuity relation
linking for the heat density $u$ with the heat flux $v$, by means of the identity
\begin{equation}
\label{continuity}
	\prt u + \partial_x v = 0.
\end{equation}
Such equation can be considered as a localized version of the global balance
\begin{equation*}
	\frac{d}{dt}\int_{C} u(t,x)\,dx+v(b)-v(a)=0,
\end{equation*}
where $C=(a,b)$ is an arbitrarily chosen control interval and $dx$ describes the length element.
Equation~\eqref{continuity} has to be coupled with a second equation relating again density $u$ and flux $v$.

\subsection{Parabolic diffusion modeling and traveling waves}

Among others, the most common choice is the {\it Fourier's law}, which is considered a good description of heat conduction,
\begin{equation}
\label{fourier}
	v=-\mu\,\partial_x u
\end{equation}
for some non-negative proportionality parameter $\mu$.
The same equation is also called {\it Fick's law} when considered in biomathematical settings, {\it Ohm's law} in electromagnetism,
{\it Darcy's law} in porous media.
In general, the coefficient $\mu$ may depend on space and time (in case of heterogeneous media) and also on the density
variable itself $u$ (and/or on its derivatives).
Here, we concentrate on the easiest case where $\mu$ is a strictly positive constant.

The coupling of~\eqref{continuity} with~\eqref{fourier} gives raise to the standard {\it parabolic diffusion equation}
\begin{equation}
\label{parabolicdiffusion}
	\prt u=\mu\,\prxx u
\end{equation}
which can be considered as a reliable description of many diffusive behaviors, such as heat conduction.
The same equation can be obtained as an appropriate limit of a brownian random walk.

Adding a reactive term $f$, which may, at first instance, depends only on the state variable $u$, consists in 
modifying the continuity equation~\eqref{continuity} into a balance law with the form
\begin{equation}
\label{balance}
	\prt u + \partial_x v=f(u).
\end{equation}
Then, coupling with the Fourier's law~\eqref{fourier}, we end up with the standard scalar parabolic
reaction--diffusion equation
\begin{equation}
	\label{reactiondiffusion}
	\prt u =\mu\,\prxx u+f(u).
\end{equation}
Two basic example of nonlinear smooth functions $f$ are usually considered
\begin{itemize}
\item[\bf i.] {\sl monostable type:} the function $f$ is strictly positive in some fixed interval, say $(U_0,U_1)$ for some $U_0<U_1$,
negative in $(-\infty,U_0)\cup(U_1,+\infty)$, and with simple zeros, i.e. $f'(U_1)<0<f'(U_0)$; 
\item[\bf ii.] {\sl bistable type.} the function $f$ is strictly positive in some fixed interval $(-\infty,U_0)\cup(U_\alpha,U_1)$
for some $U_0<U_\alpha<U_1$, negative in $(U_0,U_\alpha)\cup(U_1,+\infty)$, and with simple zeros, i.e. $f'(U_0), f'(U_1)$
strictly negative and $f'(U_\alpha)$ strictly positive.
\end{itemize}

The former case, whose prototype is $f(u)\propto u(1-u)$, corresponds to a logistic-type reaction term and it is usually
referred to as {\it Fisher--KPP equation} (using the initials of the names Kolmogorov, Petrovskii and Piscounov);
the latter, roughly given by the third order polynomial $f(u)\propto u(u-\alpha)(1-u)$ with $\alpha\in(0,1)$,
is referred to the presence of an {\it Allee-type effect} (see~\cite{CourBereGasc08}), and it is called {\it Allen--Cahn equation}
(sometimes, also named {\it Nagumo} and/or {\it Ginzburg--Landau equation}).

In both cases, the equations support existence of {\it traveling wave solutions}, namely functions with the form
$u(t,x):=\phi(\xi)$ with $\xi:=x-ct$.
Hence, the {\it profile of the wave}  $\phi$ is such that
\begin{equation*}
	\mu\,\phi''+c\,\phi'+f(\phi)=0,
\end{equation*}
for some  {\it speed} $c\in\mathbb{R}$.
Due to the fact that equation~\eqref{reactiondiffusion} is autonomous, the profile is determined up to a space translation.

In addition, traveling waves are called \\
{\bf i.} {\it traveling pulses}, if they are homoclinic orbits connecting one equilibrium with itself, that is
\begin{equation*}
	\phi_0:=\lim_{\xi\to+\infty} \phi(\xi),
\end{equation*}
for some non-constant wave profile $\phi$;\\
{\bf ii.}  {\it traveling fronts} (or {\it propagating fronts}), if they are heteroclinic orbits connecting two distinct equilibria,
that is
\begin{equation*}
	\phi_\pm:=\lim_{\xi\to+\infty} \phi(\xi).
\end{equation*}
To fix ideas, let us concentrate on the case $\phi_+$ stable.

Monotonicity of the front is a necessary condition for stability.
In fact, when dealing with partial differential equations for which a maximum principle holds, such as for the scalar
parabolic case~\eqref{reactiondiffusion}, the first eigenfunction has one sign.
Thus the first order derivative with respect to the variable $\xi$, which can be verified is an eigenfunction of the linearized
operator at the wave itself relative to the eigenvalue $\lambda=0$, is the first eigenfunction since it has one sign.
Therefore, when the maximum principle holds, all monotone waves, in case of existence, are (weakly) stable.
Analogously, non-monotone waves, again in case of existence, are unstable.

In term of existence of traveling waves, there is a significant difference between the two cases
(Fisher--KPP and Allen--Cahn), consequence of the different nature of stability of the critical points of the
associated ODE for the traveling wave profile.
Specifically, in the case of the Fisher--KPP equation, the heteroclinic orbit is a saddle/node connection;
while, in the case of the Allen--Cahn equation, it is a saddle/saddle connection.
This translates into the fact that, for the Fisher--KPP equation, there exists a (strictly negative) maximal speed $c_0$
such that traveling wave solutions exists if and only if $c\leq c_0$ (remember that we have chosen $\phi_+$ stable).
On the contrary, for the Allen--Cahn equation there exists
a unique value of the speed $c_\ast$ which corresponds to a traveling profile $\phi_\ast$.

For the Allen--Cahn equation, an explicit formula for both the profile $\phi$ and the speed $c$ can be found
in the specific case of the third order polynomial case $f(u)=\kappa\,u(u-\alpha)(1-u)$.
In this case, the equation for the traveling wave solutions can be rewritten as
\begin{equation}
\label{someq}
	\mu\,\phi''+c\,\phi'+\kappa\,\phi(\phi-\alpha)(1-\phi)=0,
\end{equation}
and thus, considering the new variable $\phi'=-A\phi(1-\phi)$ with $A>0$ to be determined, since
\begin{equation*}
	\phi''=\frac{d\phi'}{d\phi}\phi'=-A(1-2\phi)\phi'
\end{equation*}
equation~\eqref{someq} reduces to
\begin{equation*}
	\mu\,A^2(1-2\phi)+c\,A+\kappa(\phi-\alpha)=0.
\end{equation*}
Such relation can be further rewritten as a first order polynomial in $\phi$
\begin{equation*}
	(\kappa-2\mu A^2)\phi+\mu A^2+cA-\kappa\alpha=0.
\end{equation*}
In order to satisfy the identity, we need to impose the conditions
\begin{equation*}
	A=-\sqrt{\frac{\kappa}{2\mu}},\qquad c=c_\ast:=\sqrt{2\mu\kappa}\left(\frac{1}{2}-\alpha\right),
\end{equation*}
so that the unique traveling front for the Allen--Cahn equation has speed $c_\ast$
and profile $\phi$ given by the solution to
\begin{equation*}
	\phi'=-\sqrt{\frac{\kappa}{2\mu}}\phi(1-\phi)
		=-\sqrt{\frac{\kappa}{2\mu}}\phi+\sqrt{\frac{\kappa}{2\mu}}\phi^2,
\end{equation*}
which has an explicit solution given by
\begin{equation}
\label{explicitfront}
	\phi(\xi)=\frac{1}{1+e^{\sqrt{\frac{\kappa}{2\mu}}(\xi-\xi_0)}}
		=\frac{1}{2}\left\{1-\tanh(C_{\kappa,\mu}\xi)\right\},
\end{equation}
where $C_{\kappa,\mu}=\sqrt{{\kappa}/{8\mu}}$.

\subsection{Extended models}

While both the continuity equation~\eqref{continuity} and the balance law~\eqref{balance} can be considered reliable
in general contexts, the Fourier law~\eqref{fourier} should be regarded as a single possible choice among many others.
Using the same words of Onsager (cf. \cite{Onsag31a}), {\it Fourier's law is only an approximate description of the process of conduction,
neglecting the time needed for acceleration of the heat flow; for practical purposes the time-lag can be neglected in all cases
of heat conduction that are likely to be studied}.
Nevertheless, in many applications, considering extensions of the Fourier's law is required.
The first possible modification is the so-called {\it Maxwell--Cattaneo law} (or {\it Maxwell--Cattaneo--Vernotte law})
\begin{equation}
\label{maxwellcattaneo}
	\tau\prt v+v=-\mu\,\partial_x u,
\end{equation}
where $\tau>0$ is a relaxation parameter describing the time needed by the the flux $v$ to alignate
with the (negative) gradient of the density unknown $u$.
Different alternative to the Fourier's law could be considered.
Among others, let us quote here the so-called {\it Guyer--Krumhansl's law}.
In the one-dimensional setting, this consists in adding a further term
at the righthand side, namely
\begin{equation}
\label{guyerkrumhansl}
	\tau\prt v+v=-\mu\,\partial_x u+\nu\,\prxx v
\end{equation}
where $\nu>0$ is related to the mean free path of (heat) carriers. 
Both Maxwell--Cattaneo's and Guyer--Krumhansl's law can be considered as a way for incorporating into the diffusion
modeling some physical terms in the framework of Extended Irreversible Thermodynamics~\cite{CimmJouRuggVan14}.
In such a context, appropriate modification of the entropy law has to be taken into account for each one of the
corresponding modified flux laws.

Coupling~\eqref{maxwellcattaneo} with~\eqref{continuity} give raise to the classical {\it telegraph equation}
\begin{equation}
\label{telegraph}
	\tau\,\prtt u+\prt u=\mu\,\prxx u.
\end{equation}
The principal part of equation~\eqref{telegraph} coincides with the one of the wave equation, and the equation is
thus of hyperbolic type.
Therefore, for $\tau$ sufficiently small, this new equation amends a number of drawbacks inherent in~\eqref{parabolicdiffusion}
such as {\it infinite speed of propagation}, {\it ill-posedness of boundary value problems} and {\it lack of inertia}.
Here, we take into particular consideration the amendment of the latter drawback.

Similarly, coupling~\eqref{guyerkrumhansl} with~\eqref{continuity} furnishes the third order equation
\begin{equation}
\label{pseudopar}
	\tau\prtt u+\prt u=(\mu+\nu\,\prt)\prxx u,
\end{equation}
which is usually classified as a pseudo-parabolic regularization of the standard telegraph equation,
that is formally obtained in the singular limit $\nu\to 0^+$.

The variable $v$ can be eliminated from the coupled system given by the balance law~\eqref{balance}
and the Maxwell--Cattaneo equation~\eqref{maxwellcattaneo} by using the so-called \textit{Kac's trick} (see~\cite{HaMu01,Kac74}),
consisting in differentiating equation~\eqref{balance} with respect to time $t$ and the relation~\eqref{maxwellcattaneo}
with respect to space $x$ and merging them together, giving raise to the \textit{one-field equation}
\begin{equation}
\label{onefield}
	\tau\prtt u+\bigl(1-\tau f'(u)\bigr)\prt u-\mu\,\partial_{xx} u=f(u).
\end{equation}
Let us stress that the specific form for the hyperbolic reaction-diffusion equation~\eqref{onefield} depends only on the
coupling of the balance law~\eqref{balance} with the Maxwell--Cattaneo's law~\eqref{maxwellcattaneo} and not on the 
specific dependency of $f$ with respect to $u$.
In particular, the same form holds for both monostable and bistable cases.

A similar, but more complicated, equation can be in principle obtained coupling with the Guyer--Krumhansl's law,
namely
\begin{equation}
\label{guyerkrumRD}
	\tau\prtt u+\bigl(1-\tau f'(u)\bigr)\prt u=\partial_{xx}\bigl(\mu\,u-\nu f(u)+\nu\,\prt u\bigr)+f(u).
\end{equation}
which is an additional alternative pseudo-parabolic variation of~\eqref{reactiondiffusion}.

In all of the three models above presented, it is possible to introduce a convenient
rescaling of the dependent variables.
To start with, let us consider the standard reaction-diffusion equation~\eqref{reactiondiffusion}.
Next, let us introduce a rescaled variable $\tilde u$ in the form defined by
\begin{equation*}
	\tilde u=\frac{u-U_0}{U-U_0}.
\end{equation*}
for some significant value $U$.
A reasonable choice could be $U=U_1$ so that $f(U_1)=0$.
Plugging into~\eqref{reactiondiffusion}, we obtain an equation for $\tilde u$
\begin{equation*}
	\prt \tilde u =\mu\,\prxx \tilde u+\tilde f(\tilde u)
\end{equation*}
where
\begin{equation*}
	\tilde f(\tilde u):=\frac{f\bigl(U_0+(U_1-U_0)\tilde u\bigr)}{U_1-U_0},
\end{equation*}
with the advantage of having $\tilde f(1)=f(U_1)/(U_1-U_0)=0$.

Similarly, since both Maxwell--Cattaneo~\eqref{maxwellcattaneo} and Guyer--Krumhansl relations~\eqref{guyerkrumhansl}
are linear in both $u$ and $v$, considering the same scaling for $u$ and $v$
\begin{equation*}
	\tilde u=\frac{u-U_0}{U-U_0},\qquad \tilde v:=\frac{v}{U-U_0},
\end{equation*}
gives an analogous reduction to the corresponding one-field equation.
As an example, in the case of Allen--Cahn equation with relaxation, we obtain
\begin{equation*}
	\tau\prtt \tilde u+\bigl(1-\tau \tilde f'(\tilde u)\bigr)\prt \tilde u-\mu\,\partial_{xx} \tilde u=\tilde f(\tilde u),
\end{equation*}
with the same definition of $\tilde f$ reported above.
In particular, the assumption $f(1)=0$ is not restrictive.

A comprehensive theory of traveling waves for the Allen-Cahn model with relaxation is presented in~\cite{LMPS}, and further extension to the case of the Guyer-Krumhansl variation is in progress.

\subsection{Diagonalization and kinetic representation}

From now on, we will focus on the case of Allen--Cahn equation with relaxation, 
that is the semilinear hyperbolic system
\begin{equation}
\label{mainsystem}
	\prt u + \prx v = f(u)\,, \qquad \prt v + \frac{\mu}{\tau}\,\prx u = -\frac1{\tau}\,v\,,
\end{equation}
for $t\in \R^{+}$, $x \in \R$, relaxation parameter $\tau > 0$ and viscosity $\mu >0\,$, with the assumption that $f$ is of bistable type with $U_0=0$, $U_\alpha \in(0,1)$ and $U_1=1$ (refer to Section~\ref{sec:physics}). Specifically, we are interested in studying numerically the dynamics of solutions to~\eqref{mainsystem} for $f(u)=\kappa\,u(u-\alpha)(1-u)$, $\kappa>0$ and $\alpha\in(0,1)$. The corresponding Cauchy problem is determined by the initial conditions
\begin{equation}
\label{mainitialdata}
	u(0,x) = u_0(x)\,, \qquad v(0,x) = v_0(x)\,,
\end{equation}
whereas the initial conditions for~\eqref{onefield} should be assigned by deducing them from~\eqref{mainitialdata} through system~\eqref{mainsystem} as
\begin{equation*}
	u(0,x) = u_0(x)\,, \qquad \prt u(0,x) = f(u_0(x))-v_0^{\prime}(x)\,.
\end{equation*}
Setting $W=(u,v)$, together with $\flux(W)=\left(v,\frac{\mu}{\tau} u\right)$ and $\spgw(W)=\left(f(u),-\frac1{\tau} v\right)$ in~\eqref{mainsystem}, we recognize the following hyperbolic system of balance laws
\begin{equation*}
	\prt W + \prx \flux(W) = \spgw(W)\,,
\end{equation*}
where the Jacobian of the flux $\flux$ is given by the $2\! \times \!2$ constant coefficients matrix
\begin{equation*}
	\flux^{\prime} = \begin{pmatrix} 0 & 1 \\ \mu / \tau & 0\end{pmatrix},
\end{equation*}
thus leading to the nonconservative form
\begin{equation*}
	\prt W + \flux^{\prime} \prx W = \spgw(W)\,.
\end{equation*}
This system can be directly diagonalized for numerical purposes, with eigenvalues $\lambda_{\pm}=\pm \sqrt{\mu / \tau}$ and diagonalization matrix $\diag$, with its inverse $\diag^{-1}$, given by
\begin{equation*}
	\diag = \begin{pmatrix} 1 & 1 \\ -\sqrt{\mu / \tau} & \sqrt{\mu / \tau}\end{pmatrix}, \qquad 
	\diag^{-1} = \dfrac12 \!\begin{pmatrix} 1 & -\sqrt{\tau / \mu} \\ 1 & \sqrt{\tau / \mu}\end{pmatrix},
\end{equation*}
so that $\diag^{-1} \flux^{\prime}\,\diag = {\rm diag}\left(\lambda_{-},\lambda_{+}\right)$. Therefore, the diagonal variables $Z=\diag^{-1} W$, corresponding to the \textit{Riemann invariants} for the homogeneous part of~\eqref{mainsystem}, namely
\begin{equation*}
	\prt u + \prx v = 0\,, \qquad \prt v + \frac{\mu}{\tau}\,\prx u = 0\,,
\end{equation*}
they have components
\begin{equation*}
	z_{-} = \dfrac12 \left( u - \sqrt{\frac{\tau}{\mu}}\,v \right), \qquad z_{+} = \dfrac12 \left( u + \sqrt{\frac{\tau}{\mu}}\,v \right),
\end{equation*}
so that
\begin{equation}
\label{antidiag}
	u = z_{-} + z_{+}\,, \qquad v = \sqrt{\frac{\mu}{\tau}} \left(z_{+}-z_{-}\right).
\end{equation}
The source term is transformed into
\begin{equation*}
	\diag^{-1} \spgw(W) =  \dfrac12 \!\begin{pmatrix} f(u)+\frac1{\sqrt{\tau \mu}}\,v \\ f(u)-\frac1{\sqrt{\tau \mu}}\,v \end{pmatrix},
\end{equation*}
that is
\begin{equation*}
	\diag^{-1} \spgw(\diag Z) = \dfrac12 \!\begin{pmatrix} f\!\left(z_{-}+z_{+}\right)+\frac1{\tau}\left(z_{+}-z_{-}\right) \\ f\!\left(z_{-}+z_{+}\right)-\frac1{\tau} \left(z_{+}-z_{-}\right) \end{pmatrix}.
\end{equation*}
Finally, for $\varrho=\sqrt{\mu / \tau}$, the diagonal system reads
\begin{equation}
\label{diagsystem}
\left\{ \begin{aligned}
	\prt z_{-} - \varrho\,\prx z_{-} = \frac12\,f\!\left(z_{-}+z_{+}\right) + \frac1{2\tau}\left(z_{+}-z_{-}\right)\\
	\prt z_{+} + \varrho\,\prx z_{+} = \frac12\,f\!\left(z_{-}+z_{+}\right) - \frac1{2\tau}\left(z_{+}-z_{-}\right)
\end{aligned} \right.
\end{equation}
meaning that the diagonal variables satisfy the so-called weakly coupled semilinear \textit{Goldstein--Taylor model} of diffusion equations. Such system admits an important physical interpretation, since it can be interpreted as the reactive version of the hyperbolic \textit{Goldstein--Kac model}~\cite{Kac74} for the (easiest possible) correlated random walk. In view of its numerical approximation, this representation is intrinsically \textit{upwind} in the sense that $z_{-}$ represents the contribution to the density $u$ of the particles moving to the left with negative velocity $-\varrho\,$, while $z_{+}$ corresponds to the particles moving to the right with positive velocity $\varrho\,$, according to the uniform jump process with equally distributed transition probability.


\section{Formulation of the numerical method}
\label{sec:numerics}

We perform \textit{finite volume schemes} because of the possible implementation for models with low regularity of the solutions, so that an integral formulation is suitable. Moreover, nonuniform discretizations of the physical space are specially required, taking into account the typical inhomogeneity of the dynamics over different regions. This is important as well for computational issues, when nonuniform time-grids are used for improving the CPU performance.

\subsection{First order scheme and nonuniform grids}

\begin{figure*}[tb]
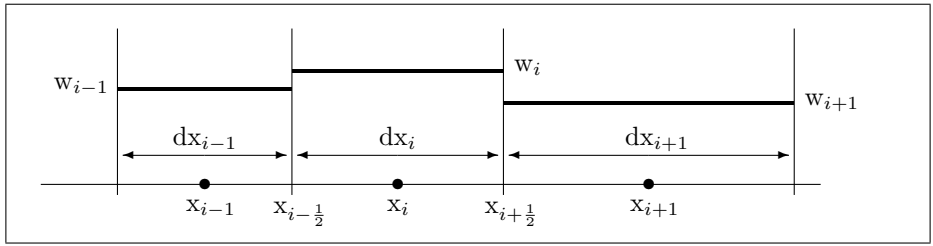

\blankbox{.97\columnwidth}{7.5pc
\put(42,18){\line(0,1){63}}
\put(13,22){\line(1,0){295}}
\put(75,22){\circle*{4}} \put(68,12){$\xx_{i-1}$}
\put(75,33){\vector(-1,0){31}}
\put(75,33){\vector(1,0){31}}
\put(63,37){$\dxim$} \put(101,11){$\xim$}
\put(108,18){\line(0,1){63}}
\put(148,22){\circle*{4}} \put(144,12){$\xx_i$}
\put(148,33){\vector(-1,0){38}}
\put(148,33){\vector(1,0){38}}
\put(141,37){$\dxi$}
\put(188,18){\line(0,1){63}}
\put(243,22){\circle*{4}} \put(236,12){$\xx_{i+1}$}
\put(298,18){\line(0,1){63}}
\put(243,33){\vector(-1,0){53}}
\put(243,33){\vector(1,0){53}}
\put(234,37){$\dxip$} \put(181,11){$\xip$}
\linethickness{0.4mm}
\put(42,58){\line(1,0){66}} \put(18,58){${\rm w}_{i-1}$}
\put(108,65){\line(1,0){80}} \put(192,65){${\rm w}_i$}
\put(188,53){\line(1,0){110}} \put(302,52){${\rm w}_{i+1}$}}
\caption{piecewise constant reconstruction on nonuniform mesh/grid \eqref{interfacial}}
\label{fig:firstmesh}
\end{figure*}

We set up a nonuniform mesh on the one-dimensional space (see Figure~\ref{fig:firstmesh}) and we denote by $C_i\!=\![\xim,\xip)$ the finite volume (cell) centered at point $\xx_i\!=\!\tfrac{1}{2}(\xim+\xip),\,i\!\in\!\Z\,$, where $\xim$ and $\xip$ are the cell's interfaces and $\dxi\!=\!\text{length}(C_i)$, therefore the characteristic space-step is given by ${\rm dx}\!=\!\sup_{i\in\Z} \dxi\,$. We build a piecewise constant approximation of any (sufficiently smooth) function by means of its \textit{integral cell-averages}, namely
\begin{equation}
\label{effenum}
{\rm w}_i = \frac1{\dxi} \int_{C_i}\!w(x)\,dx \approx w(\xx_i) + {\mathcal O}({\rm dx}^2)\,,
\end{equation}
because of the symmetric integral $\int_{C_i}(x-\xx_i)\,dx=0$ due to the cell-centered structure of the grid, that converges uniformly to~$w(x)$ as ${\rm dx} \rightarrow 0\,$. Moreover, a straightforward computation leads to the approximation
\begin{equation}
\label{interfacial}
	{\rm w}_{i+1} - {\rm w}_i \,=\, w^{\prime}(\xx_i)\!\left( \!\frac{\dxip}2 + \frac{\dxi}2 \!\right) 
	+ {\mathcal O}({\rm dx}^2)\,,
\end{equation}
that is defined over an interfacial interval $[\xx_i,\xx_{i+1}]$ and, for example, it reproduces the correct \textit{upwind interfacial quadrature} for the advection with negative speed if we observe that
\begin{equation}
\label{effederiv}
	\frac1{\dxi} \int_{C_i}\!w^{\prime}(x)\,dx = \frac1{\dxi} \left( w(\xip) - w(\xim) \right).
\end{equation}

In that framework, a \textit{semi-discrete finite volume scheme} applied to the system~\eqref{diagsystem} produces a numerical solution in the form of a (discrete valued) vector whose in-cell values are interpreted as approximations of the cell-averages, i.e.
\begin{equation}
\label{averages}
\rd_i(t) \approx \frac1{\dxi} \int_{C_i}\!z_{-}(t,x)\,dx \,, \quad \sd_i(t) \approx \frac1{\dxi} \int_{C_i}\!z_{+}(t,x)\,dx \,,
\end{equation}
and which satisfy the upwind three-points scheme
\begin{equation}
\label{semidiag}
\begin{split}
	\frac{d \rd_i}{dt} &= \frac{\varrho}{\dxi} \left(\rd_{i+1} - \rd_i \right) + \frac12\,f\!\left(\rd_i+\sd_i\right) 
		+ \frac1{2\tau} \left(\sd_i - \rd_i\right)\\
	\frac{d \sd_i}{dt} &= -\frac{\varrho}{\dxi} \left(\sd_i - \sd_{i-1} \right) + \frac12\,f\!\left(\rd_i+\sd_i\right) 
		- \frac1{2\tau} \left(\sd_i - \rd_i\right)
\end{split}
\end{equation}
when considering~\eqref{effederiv} for the diagonal variables in~\eqref{diagsystem} which are advected with constant speed. By setting $\uu_i = \rd_i+\sd_i$ and $\vv_i = \varrho \left(\sd_i - \rd_i \right)$ according to~\eqref{antidiag}, and recalling that $\varrho=\sqrt{\mu/\tau}\,$, we obtain through a straightforward computation a semi-discrete version of~\eqref{mainsystem} that is
\begin{equation}
\label{semimain}
\begin{split}
	\frac{d \uu_i}{dt} &= - \frac{\vv_{i+1} - \vv_{i-1}}{2 \dxi} + f(\uu_i) 
		+ \frac12 \varrho\,\dxi\,\frac{\uu_{i+1}-2\uu_i+\uu_{i-1}}{\dxi^2}\\
	\frac{d \vv_i}{dt} &= - \varrho^2 \frac{\uu_{i+1} - \uu_{i-1}}{2 \dxi} - \frac1{\tau}\,\vv_i 
		+ \frac12 \varrho\,\dxi\,\frac{\vv_{i+1}-2\vv_i+\vv_{i-1}}{\dxi^2}
\end{split}
\end{equation}
with initial data corresponding to~\eqref{mainitialdata} by means of an approximate condition
\begin{equation*}
\uu_i(0) = \frac1{\dxi} \int_{C_i}\!u_0(x)\,dx \,, \quad \vv_i(0) = \frac1{\dxi} \int_{C_i}\!v_0(x)\,dx \,, \qquad i \in \Z\,.
\end{equation*}

\noindent It is worthwhile noticing that, in case of uniform grids, i.e. $\dxi={\rm dx}\,$, for any $i\in \Z\,$, a standard Taylor's expansion from~\eqref{effenum}-\eqref{interfacial} leads to show that~\eqref{semimain} formally corresponds to
\begin{equation*}
	\prt u + \prx v = f(u) + \frac12 \varrho\,{\rm dx}\,\prxx u\,, \qquad \prt v + \varrho^2 \prx u = - \frac1{\tau}\,v + \frac12 \varrho\,{\rm dx}\,\prxx v\,,
\end{equation*}
so that the scheme is consistent in the usual sense of the \textit{modified equation}~\cite{LV}, although we expect the appearance of a numerical viscosity with strength measured through the physical and numerical parameters $\varrho$ and ${\rm dx}\,$.

However, the utilization of unstructured spatial grids is required for problems incorporating composite geometries, also in view of the recent theoretical advances on adaptive techniques for mesh refinement in the resolution of multi-scale complex systems. For the case of a nonuniform mesh, the approximation~\eqref{interfacial} seems to reveal a lack of consistency of the numerical scheme~\eqref{semimain} with the underlying continuous equations, as the space-step $\dxi$ could be very different from the length of an interfacial interval $\bigl|\xx_{i+1}-\xx_i\bigr|=\tfrac{1}{2}\dxi +\tfrac{1}{2}{\dxip}$. Nevertheless, the issue of an error analysis with optimal rates can be pursued, by virtue of the results concerning the \textit{supra-convergence phenomenon} for numerical approximation of hyperbolic conservation laws. In fact, despite a deterioration of the pointwise consistency is observed in consequence of the non-uniformity of the mesh, the formal accuracy is actually maintained as the global error behaves better than the (local) truncation error would indicate. This property of enhancement of the numerical error has been widely explored, and the question of (finite volume) upwind schemes for conservation laws and balance equations is addressed in~\cite{BGP}, \cite{KS} and~\cite{CHS}, with proof of convergence at optimal rates for smooth solutions.

\subsection{Time discretization}

We introduce a variable time-step $\dtn\!=\!{\rm t}_{n+1}\!-\!{\rm t}_n$, $n\!\in\!\N$, and we set ${\rm dt}\!=\!\sup_{n\in\N} \dtn\,$, therefore we have to consider a {\sl CFL-condition}~\cite{LV} on the ratio $\frac{\dtn}{\dxi}$ for the numerical stability. We discretize the time operator in~\eqref{semidiag} by means of a mixed explicit-implicit approach, as follows
\begin{equation*}
\begin{split}
	\frac{\rd_i^{n+1}-\rd_i^{n}}{\dtn} &= \frac{\varrho}{\dxi} \bigl(\rd_{i+1}^{n+1}-\rd_i^{n+1}\bigr) 
	+ \frac12\,f(\rd_i^{n}+\sd_i^{n}) + \frac1{2\tau} \bigl(\sd_i^{n+1}-\rd_i^{n+1}\bigr)\\
	\frac{\sd_i^{n+1}-\sd_i^{n}}{\dtn} &= -\frac{\varrho}{\dxi} \bigl(\sd_i^{n+1}-\sd_{i-1}^{n+1}\bigr) 
	+ \frac12\,f(\rd_i^{n}+\sd_i^{n}) - \frac1{2\tau} \bigl(\sd_i^{n+1}-\rd_i^{n+1}\bigr)
\end{split}
\end{equation*}
Fully implicit schemes have also been tested with no significant advantage in the quality of the approximation, but with a significant increase of the computational time.

At this point, an important simplification in terms of the actual implementation of the above algorithm arises if considering uniform time and space stepping, i.e. $\dtn={\rm dt}\,$, for any $n\in \N\,$ and $\dxi={\rm dx}\,$, for any $i\in \Z\,$. Indeed, by setting
\begin{equation*}
	\alpha = \varrho\frac{\rm dt}{\rm dx}\,, \qquad \beta = \frac{\rm dt}{2\tau}\,, \qquad {\rm f}_i^n = f(\rd_i^n+\sd_i^n)\,,
\end{equation*}
the above algorithm can be rewritten in compact form as
\begin{equation}
\label{vectorscheme}
	\begin{pmatrix}
	(1+\beta)\mathbb{I}-\alpha\,\mathbb{D}_{+} & - \beta\,\mathbb{I} \\
	- \beta\,\mathbb{I} & (1+\beta)\mathbb{I}+\alpha\,\mathbb{D}_{-} \\
	\end{pmatrix}
	\!\begin{pmatrix} \rd^{n+1} \\ \sd^{n+1} \end{pmatrix} = 
	\begin{pmatrix} \rd^n + \frac{\rm dt}2{\rm f}^n \\ \sd^n + \frac{\rm dt}2{\rm f}^n \end{pmatrix}
\end{equation}
where the matrices $\mathbb{I}\,$, $\mathbb{D}_{-}$ and $\mathbb{D}_{+}$ are given by
\begin{equation*}
	\mathbb{I}=(\delta_{i,j})\,, \qquad 
	\mathbb{D}_{-} = \bigl( \delta_{i,j} - \delta_{i,j+1} \bigr)\,, \quad 
	\mathbb{D}_{+} = \bigl( \delta_{i+1,j} - \delta_{i,j} \bigr)\,,
\end{equation*}
and $\delta_{i,j}$ is the standard \textit{Kronecker symbol}\,. The block-matrix in~\eqref{vectorscheme} is invertible, since its spectrum is contained in the complex half plane $\bigl\{ \lambda \in \mathbb{C}\,:\,{\rm Rel}(\lambda) \geq 1 \bigr\}$ as a consequence of the \textit{Ger\v sgorin criterion}~\cite{QSS}.\\
A direct manipulation of~\eqref{vectorscheme} gives
\begin{equation}
\label{finalscheme}
	\begin{aligned}
	\rd^{n+1} &= \bigl( \mathbb{S} - \alpha^2\,\mathbb{D}_{-} \mathbb{D}_{+} \bigr)^{-1} 
		\Bigl\{ \bigl[ (1+\beta)\mathbb{I}+\alpha\,\mathbb{D}_{-} \bigr] \rd^n + \beta\,\sd^n\\
		& \hskip4.5cm + \frac{\rm dt}2 \bigl[ (1+2\beta)\mathbb{I} + \alpha\,\mathbb{D}_{-} \bigr] {\rm f}^n \Bigr\}\\
	\sd^{n+1} &= \bigl( \mathbb{S} - \alpha^2\,\mathbb{D}_{+} \mathbb{D}_{-} \bigr)^{-1} 
		\Bigl\{ \beta\,\rd^n + \bigl[ (1+\beta)\mathbb{I} - \alpha\,\mathbb{D}_{+} \bigr] \sd^n\\
		& \hskip4.5cm + \frac{\rm dt}2 \bigl[ (1+2\beta)\mathbb{I} - \alpha\,\mathbb{D}_{+} \bigr] {\rm f}^n \Bigr\}
	\end{aligned}
\end{equation}
where $\mathbb{S}$ is the symmetric matrix
\begin{equation*}
	\mathbb{S} = (1+2\beta)\mathbb{I}+\alpha\,(1+\beta)\bigl(\mathbb{D}_{-} - \mathbb{D}_{+}\bigr)\,.
\end{equation*}

Nevertheless, one of the most important features of the models described in Section~\ref{sec:physics} is that they could produce strikingly nontrivial patterns. Therefore, the use of nonuniform meshes is somehow mandatory and hence the numerical solution often requires very long computational time, for the large amount of data to be traded in order to accurately capture the details of physical phenomena. Moreover, especially for applied scientists involved in setting up realistic experiments, the possibility of running fast comparative simulations using simple algorithms implemented into affordable processors is of a primary interest. In this context, parallel computing based on modern graphics processing units (GPUs) enjoys the advantages of a high performance system with relatively low cost, allowing for software development on general-purpose microprocessors even in personal computers. As a matter of fact, GPUs are revolutionizing scientific simulation by providing several orders of magnitude of increased computing capability inside a mass-market product, making these facilities economically attractive across subsets of industry domains~\cite{HHH,TNL,MIM,HWW}. Simple approximation schemes like~\eqref{semimain} are often acceptable even for real problems, so that proper numerical modeling becomes accessible to practitioners from various scientific fields.

\subsection{Second order scheme}

\begin{figure*}[tb]
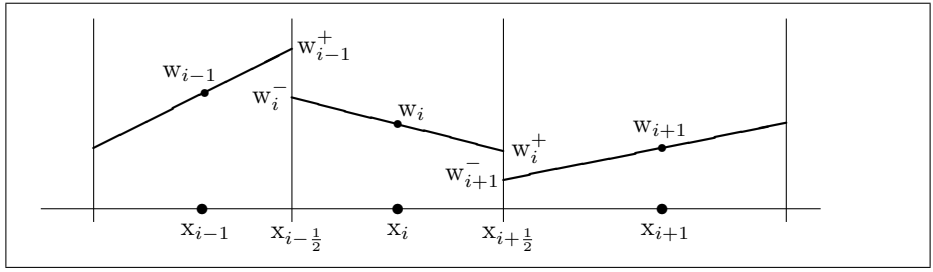

\blankbox{.97\columnwidth}{8.3pc
\put(13,22){\line(1,0){295}}
\put(33,17){\line(0,1){77}}
\put(295,17){\line(0,1){77}}
\put(74,22){\circle*{4}} \put(66,12){$\xx_{i-1}$}
\put(108,17){\line(0,1){77}} \put(100,11){$\xim$}
\put(148,22){\circle*{4}} \put(144,12){$\xx_i$}
\put(188,17){\line(0,1){77}} \put(180,11){$\xip$}
\put(248,22){\circle*{4}} \put(240,12){$\xx_{i+1}$}
\put(75,66){\circle*{3}} \put(59,72){${\rm w}_{i-1}$}
\put(148,54){\circle*{3}} \put(148,57){${\rm w}_i$}
\put(248,45){\circle*{3}} \put(237,51){${\rm w}_{i+1}$}
\thicklines
\put(33,45){\line(2,1){75}}
\put(108,64){\line(4,-1){80}}
\put(188,33){\line(5,1){107}}
\put(93,62){${\rm w}_i^{-}$}
\put(191,42){${\rm w}_i^{+}$}
\put(110,81){${\rm w}_{i-1}^{+}$}
\put(166,33){${\rm w}_{i+1}^{-}$}}
\caption{piecewise linear reconstruction on nonuniform mesh}
\label{fig:secondmesh}
\end{figure*}

The basic idea to develop second order schemes is to replace the piecewise constant reconstruction~\eqref{effenum} by piecewise linear approximations (see Figure~\ref{fig:secondmesh}), which provide more accurate values at the cell's interfaces.

On that account, based on the cell-averages, we associate to~\eqref{averages} some correct coefficients, for all $i\!\in\! \Z\,$, $x\!\in\! C_i\,$, which are given by
\begin{equation}
\label{effenum2}
	\rd_i(t,x) = \rd_i(t)+(x-\xx_i)\,\rd_i^{\prime}\,, \qquad \sd_i(t,x) = \sd_i(t)+(x-\xx_i)\,\sd_i^{\prime}\,,
\end{equation}
where $\rd^{\prime}_i$ and $\sd^{\prime}_i$ indicate the numerical derivatives, which are defined as appropriate interpolations of the discrete increments between neighboring cells, for example,
\begin{equation}
\label{slopelimiter1}
	\rd_i^{\prime} = lmtr \left\{ \frac{\rd_{i+1}-\rd_i}{\xx_{i+1}-\xx_i} \,, \frac{\rd_i-\rd_{i-1}}{\xx_i-\xx_{i-1}} \right\}, \quad i \in \Z\,.
\end{equation}
Because also higher-order reconstructions are, in general, discontinuous at the cell's interfaces, possible oscillations are suppressed by applying suitable \textit{slope-limiter} techniques (see~\cite{HO,HH} for instance).

Therefore, second order interpolations are computed from~\eqref{effenum2} to define the interfacial values at $\xim$ and $\xip$ as follows
\begin{align*}
	\rd_i^{-}(t) & = \rd_i(t) - \frac{\dxi}2\,\rd_i^{\prime}\,, \quad \rd_i^{+}(t) = \rd_i(t) + \frac{\dxi}2\,\rd_i^{\prime}\,, \\
\sd_i^{-}(t) & = \sd_i(t) - \frac{\dxi}2\,\sd_i^{\prime}\,, \quad \sd_i^{+}(t) = \vv_i(t) + \frac{\dxi}2\,\sd_i^{\prime}\,,
\end{align*}
which are then substituted inside~\eqref{semidiag} to obtain more accurate numerical jumps at the interfaces, namely
\begin{equation}
\label{semidiag2}
\begin{split}
	\frac{d \rd_i}{dt} &= \frac{\varrho}{\dxi} \left(\rd_{i+1}^{-} - \rd_i^{+} \right) 
	+ \frac12\,f\!\left(\rd_i+\sd_i\right) + \frac1{2\tau} \left(\sd_i - \rd_i\right)\\
	\frac{d \sd_i}{dt} &= -\frac{\varrho}{\dxi} \left(\sd_i^{-} - \sd_{i-1}^{+} \right) 
	+ \frac12\,f\!\left(\rd_i+\sd_i\right) - \frac1{2\tau} \left(\sd_i - \rd_i\right)
\end{split}
\end{equation}
We notice that, the equation being linear in the principal hyperbolic part, the second order scheme with \textit{flux limiter} in~\cite{BCN} is precisely of the type above, since the flux is trivially given by the conservation variables.

For the sake of simplicity, we have been considering in Section~\ref{sec:experiments} only the first order discretization in time, but it is easy recovering higher order accuracy by applying Runge-Kutta methods (refer to~\cite{GS} for an overall introduction), that appears to be essential for practical computations.


\section{Numerical simulations}
\label{sec:experiments}

We start by briefly revising some basics of the numerical results in~\cite{LMPS}, in order to assess the reliability of the numerical method presented in Section~\ref{sec:numerics} for determining the behavior of the solutions to reaction-diffusion models with relaxation introduced in Section~\ref{sec:physics}.

We use the algorithm~\eqref{finalscheme} to analyse the wave speeds $c_\ast$ of the traveling front connecting the stable states $0$ and $1$. Following~\cite{LeVYee90}, we introduce an \textit{average speed} of the numerical solution
at time ${\rm t}^n$ defined by
\begin{equation}
\label{numerAve}
  c^n = \frac1{\rm dt} \mathbf{1} \cdot (\uu^{n} - \uu^{n+1}) 
  	= \frac1{\rm dt} \sum_i ( \uu^{n}_i - \uu^{n+1}_i ),
\end{equation}
where $\mathbf{1}=(1,\dots,1)$ and recalling that $\uu^n=\rd^n+\sd^n\,, n\in \N\,$. We consider the bistable function $f(u)=u(u-\alpha)(1-u)$ with $\alpha\in(0,1)$, aiming at comparing the values for the propagation speed $c_\ast$ as obtained by means
of the shooting argument in~\cite{LMPS} and the ones given by~\eqref{numerAve}.

The solution to the Cauchy problem is selected with an increasing
datum connecting $0$ and $1\,$, and then computing $c^n$ at a time $t$ so large that stabilization
of the propagation speed for the numerical solution is reached. We have been testing three choices for the couple $(\tau, \alpha)$ for different values
of ${\rm dt}$ and ${\rm dx}$, where the range of variation of $\tau$ is chosen so that the condition
$\tau\,f'(u)<1$ is satisfied for all values of the unstable zero $\alpha$ (see Table~\ref{tab:numerspeed}). Requiring to detect the speed value with an error always less than
5\% of the effective value, we heuristically determine ${\rm dx}=2^{-3}$ and
${\rm dt}=10^{-2}$, that will be used for subsequent numerical experiments. 
For such a choice, we record in Table~\ref{tab:firstorder} the results of the first order scheme 
for various values of $\alpha$ and $\tau=1$ or $\tau=4$ (together with
the corresponding relative error) and in Table~\ref{tab:secondorder} those of a second order scheme.

\begin{table}
\centering
\caption{Relative error for the numerical velocity of the Riemann problem with jump at $\ell/2$, $\ell=25$ ($T$ final time and $N$ number of grid points):
A. $\tau=1$, $\alpha=0.9$, $c_\ast=0.5646$, $T=40$;
B. $\tau=2$, $\alpha=0.6$, $c_\ast=0.1737$, $T=30$;
C. $\tau=4$, $\alpha=0.7$, $c_\ast=0.3682$, $T=35$.}
\vskip6pt
{\begin{tabular}{@{}r|c|r|r|r|r|r@{}} 
 			&${\rm dx}$ 	&$2^0$	&$2^{-1}$	&$2^{-2}$	&$2^{-3}$	&$2^{-4}$	\\ \hline
 			&A 		&0.1664	&0.0787	&0.0325	&0.0091	&0.0018 	\\
${\rm dt}=10^{-1}$ 	&B 		&0.0383	&0.0306	&0.0241	&0.0198	&0.0175 	\\
 			&C 		&0.1527 	&0.1144	&0.0818	&0.0581	&0.0442 	\\ 
\hline
			&A 		&0.1751	&0.0876	&0.0417	&0.0186	&0.0079	\\
${\rm dt}=10^{-2}$	&B 		&0.0275	&0.0196	&0.0128	&0.0084	&0.0061	\\
			&C 		&0.1420	&0.1018	&0.0684	&0.0457	&0.0339	\\
\hline
			&A 		&0.1760 	&0.0885	&0.0427	&0.0196	&0.0089	\\
${\rm dt}=10^{-3}$	&B 		&0.0265 	&0.0184 	&0.0117	&0.0072	&0.0049	\\
			&C 		&0.1411	&0.1006	&0.0670	&0.0441	&0.0321
\end{tabular}}
\label{tab:numerspeed}
\end{table}

\begin{table}
\centering
\caption{First order in space: final average speed~\eqref{numerAve} and relative error
with respect to $c_\ast$ given in~\cite{LMPS} 
($N=400$, ${\rm dx}=0.125$, ${\rm dt}=0.01$, $\ell=25$, $T=40$)}
\vskip6pt
{\begin{tabular}{@{}r|c|c|c|c@{}} 
 			&$\alpha=0.6$	&$\alpha=0.7$	&$\alpha=0.8$	&$\alpha=0.9$	\\ \hline
$\tau=1$ 		&0.1580		&0.3096		&0.4497		&0.5751 	\\
	&{\scriptsize 0.0101}	 &{\scriptsize 0.0118}  &{\scriptsize 0.0145} &{\scriptsize 0.0186} \\
\hline
$\tau=4$		&0.2102		&0.3533		&0.4337		&0.4825	\\
	&{\scriptsize 0.0396} &{\scriptsize 0.0404} &{\scriptsize 0.0365} &{\scriptsize 0.0118}
\end{tabular}}
\label{tab:firstorder}
\end{table}

\begin{table}
\centering
\caption{Second order in space: final average speed~\eqref{numerAve} and relative error
with respect to $c_\ast$ given in~\cite{LMPS} 
($N=400$, ${\rm dx}=0.125$, ${\rm dt}=0.01$, $\ell=25$, $T=40$)}
\vskip6pt
{\begin{tabular}{@{}r|c|c|c|c@{}} 
 			&$\alpha=0.6$	&$\alpha=0.7$	&$\alpha=0.8$	&$\alpha=0.9$	\\ \hline
$\tau=1$ 		&0.1560		&0.3052		&0.4421	 	&0.5630	\\
	&{\scriptsize 0.0025}	  &{\scriptsize 0.0025}  &{\scriptsize 0.0026}  &{\scriptsize 0.0029} \\
\hline
$\tau=4$		&0.2184		&0.3672		&0.4485		&0.4885	\\
	&{\scriptsize 0.0022}  &{\scriptsize 0.0025}  &{\scriptsize 0.0034}  &{\scriptsize 0.0004}
\end{tabular}}
\label{tab:secondorder}
\end{table}

\subsection{Riemann problem as a large perturbation}

For these applications, we restrict to the first order discretization, since we are interested in considering initial data 
with sharp transitions. In such cases, higher order approximations of the derivatives typically 
introduce spurious oscillations that, even being transient and possibly cured by employing suitable \textit{slope limiters}\,, they may however lead to catastrophic
consequences because of the bistable nature of the reaction term.

The main achievement is 
that we are able to show that the actual domain of attraction of the
front is much larger than guaranteed by the nonlinear stability analysis performed in~\cite{LMPS}.
Indeed, the analytical results state that small perturbations
to the propagating fronts are dissipated, with an exponential rate. Nevertheless, we expect that 
the front possesses a larger domain of attraction (as already known for the parabolic Allen--Cahn equation~\cite{FifeMcLe77}) and, specifically, that any bounded initial data $u_0$ such that 
\begin{equation}
\label{decay}
	\limsup_{x\to-\infty} \,u_0(x) < \alpha < \liminf_{x\to+\infty} \,u_0(x)
\end{equation}
gives raise to a solution that is asymptotically convergent to some traveling
front connecting $u=0$ with $u=1$.

\begin{figure}[hbt]
\begin{center}
{\includegraphics[width=6.75cm]{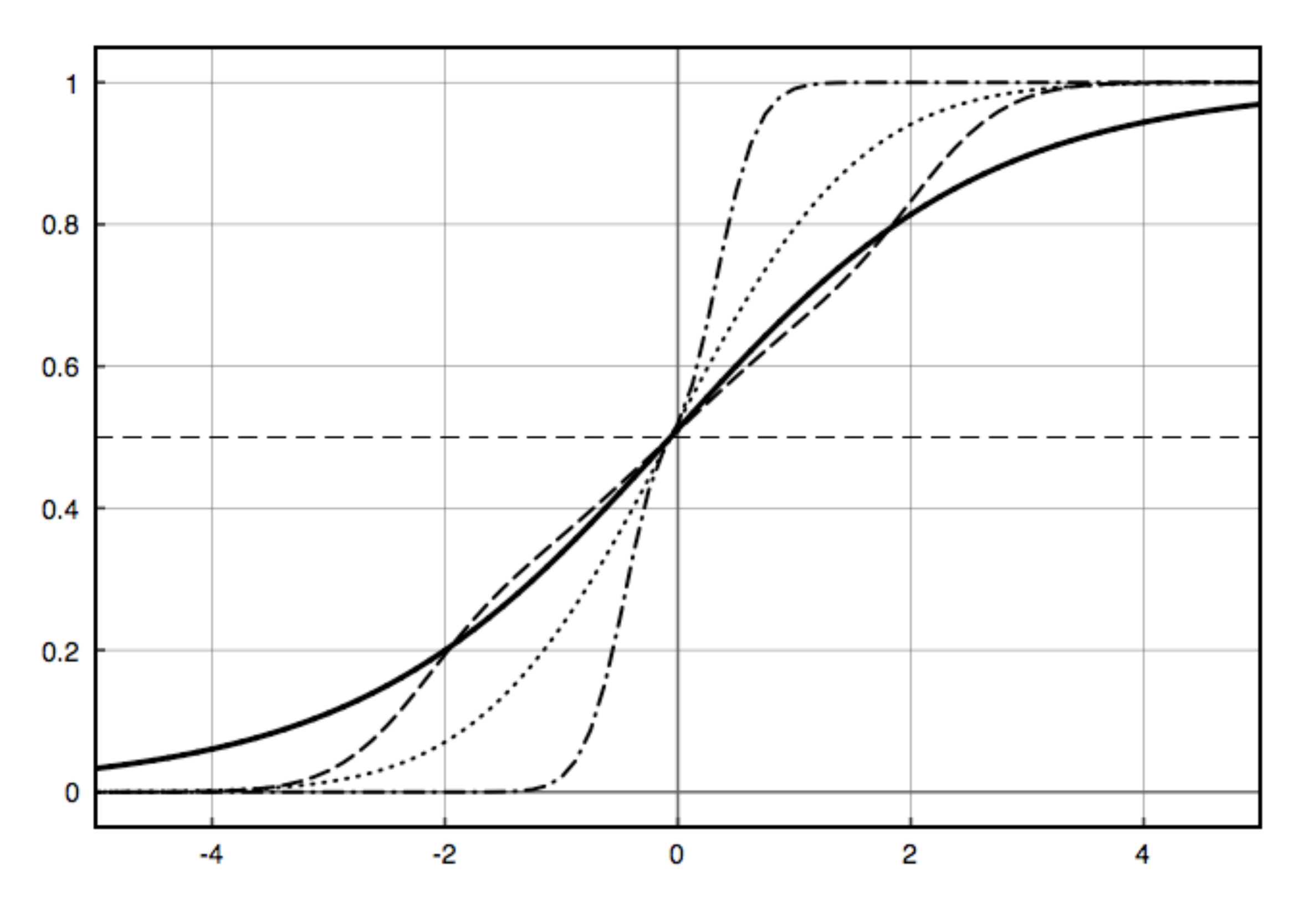}}
{\includegraphics[width=6.75cm]{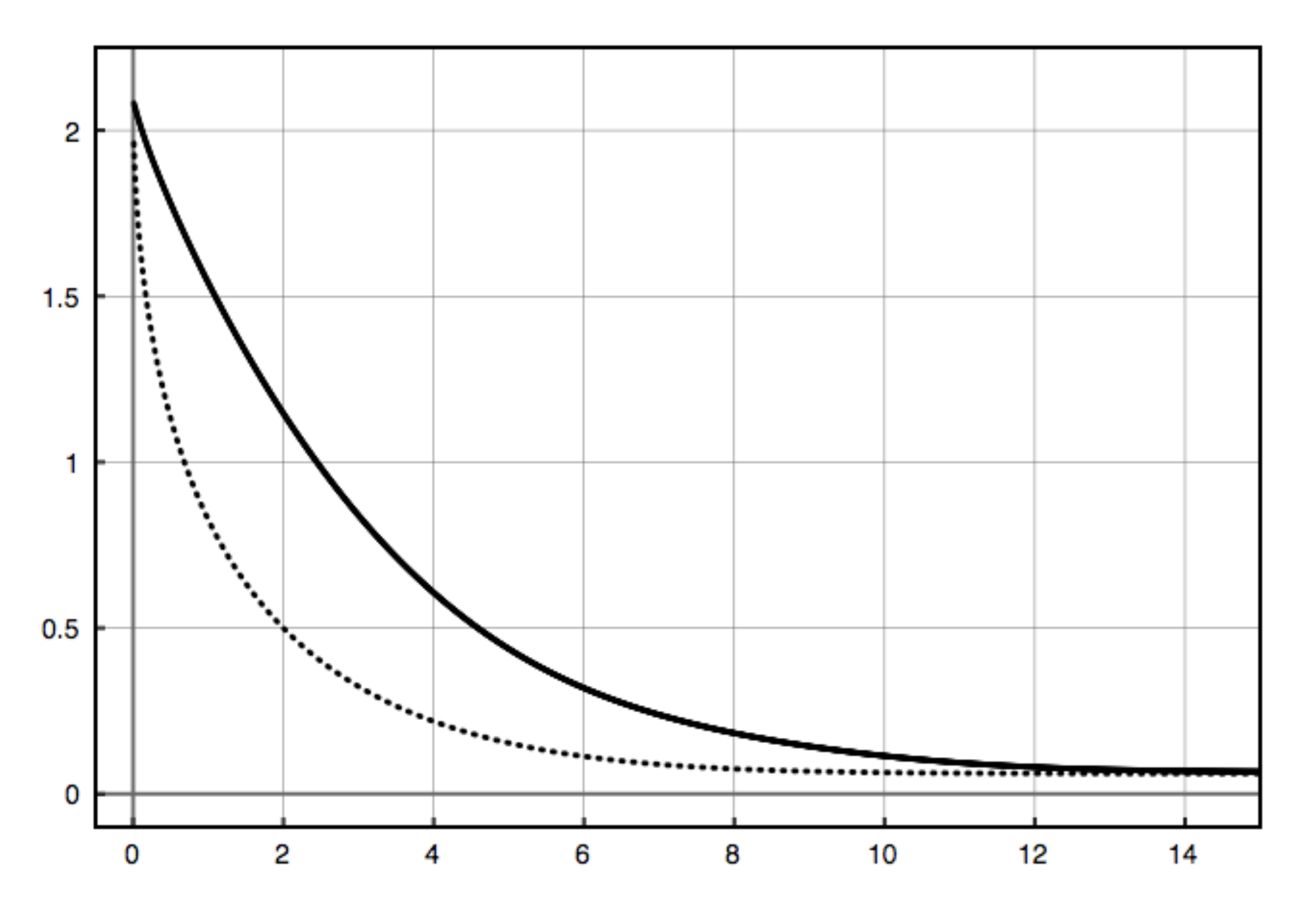}}
\end{center}
\caption{Riemann problem with initial datum $\chi_{{}_{(0,\ell)}}$ in $(-\ell,\ell)$, $\ell=25$.
Left: solution profiles zoomed in the interval $(-5,5)$ at time $t=1$ (dash-dot), $t=5$ (dash), $t=15$ (continuous), for comparison, solution to the parabolic Allen--Cahn equation at time $t=1$ (dot).
Right: Decay of the $L^2$ distance to the exact equilibrium solution for the hyperbolic
(continuous) and parabolic (dot) Allen--Cahn equations.}
\label{fig:Riemann2}
\end{figure}

To support such conjecture, we perform numerical experiments with
\begin{equation*}
	\tau=4\,, \qquad \ell=25\,, \qquad {\rm dx}=0.125\,, \qquad {\rm dt}=0.01\,.
\end{equation*}
We consider the case $\alpha=1/2$ motivated by the fact that the profile of the traveling front for the hyperbolic Allen--Cahn equation is stationary and it coincides
with the one of the corresponding original parabolic equation, explicitly given by~\eqref{explicitfront} and normalized by the condition $U(0)=1/2$. Numerical simulations confirm the decay of the solution to the equilibrium profile (see Figure~\ref{fig:Riemann2}, left). When compared with the standard Allen--Cahn equation, it appears evident that the dissipation mechanism of the hyperbolic equation is weaker with respect to the parabolic case (see Figure~\ref{fig:Riemann2}, right).

\subsection{Randomly perturbed initial data}

The genuine novelty of the numerical simulations illustrated in this section consists in suggesting that stability of the traveling waves actually goes beyond the regime $1-\tau f'(u)$ positive, that is required in the theoretical statements proven in~\cite{LMPS}.

\begin{figure}[hbt]
\begin{center}
{\includegraphics[width=6.75cm]{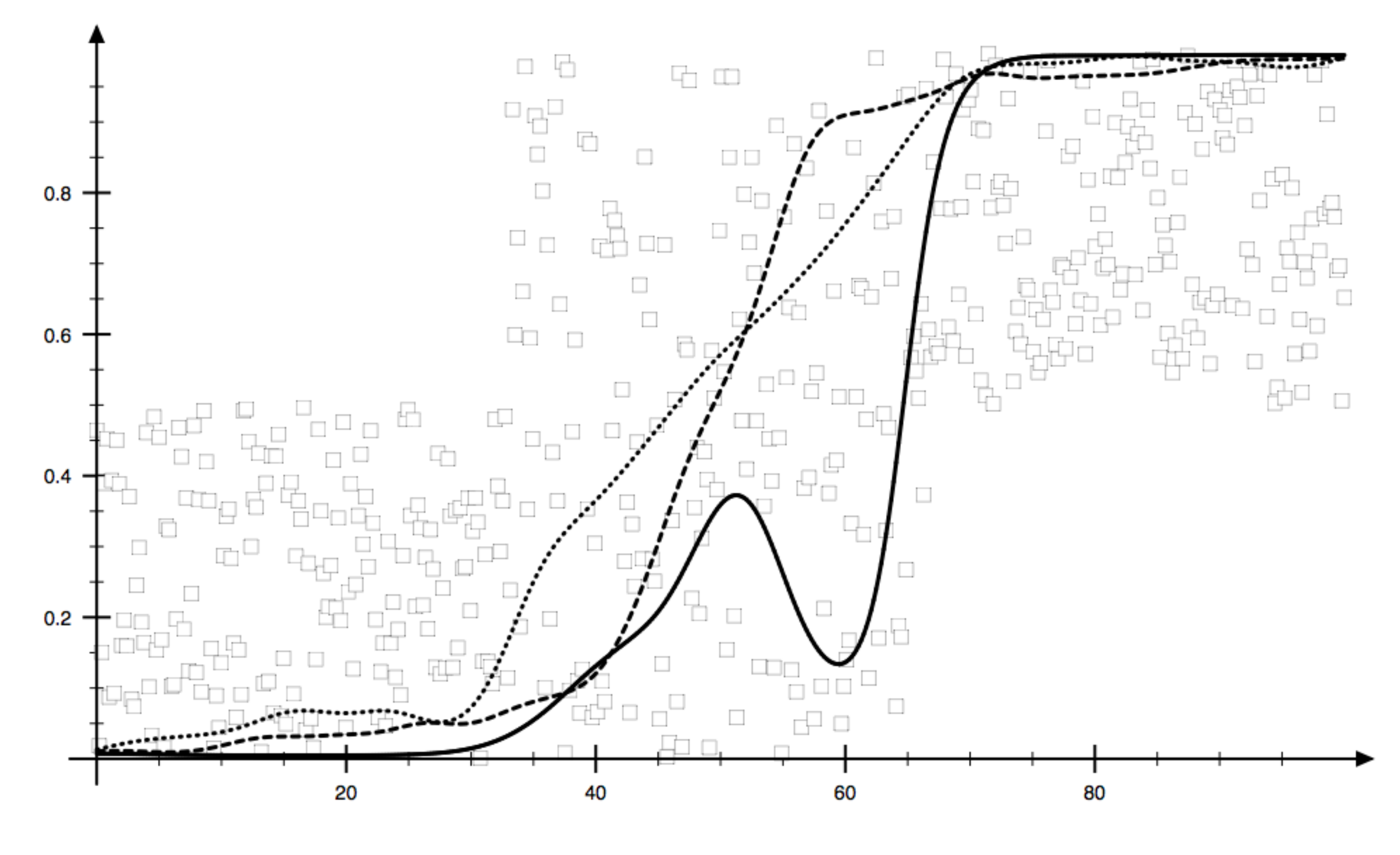}}
{\includegraphics[width=6.75cm]{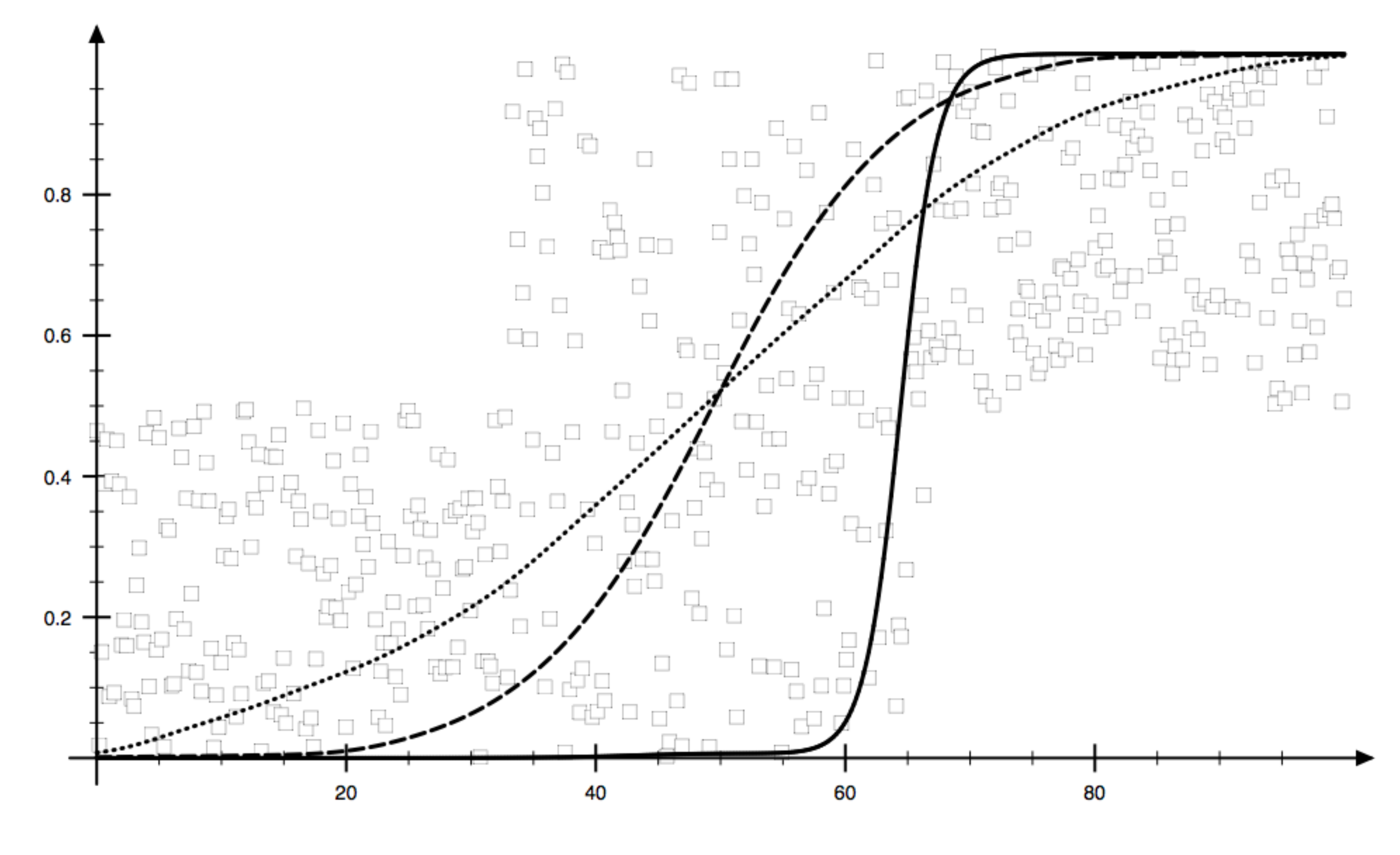}}
\end{center}
\caption{Random initial datum ($\square$).
Solution profiles for the hyperbolic Allen--Cahn equation with relaxation at time $t=10$ (left)
and time $t=20$ (right) for $\tau=1$ (continuous line), $\tau=5$ (dashed) and $\tau=10$ (dotted).
}\label{fig:Random1tau}
\end{figure}

We consider initial data that resemble only very roughly the transition from 0 to 1.
More precisely, we divide the interval $(0,\ell)$ into three parts of equal length and we choose a 
random value in each of these sub-intervals coherently with the requirement~\eqref{decay}.
We assign $u_0(x)$ to be any different random value in $(0,0.5)$ for each $x\in(0,\ell/3)$,
in $(0,1)$ for each $x\in(\ell/3,2\ell/3)$ and in $(0.5,1)$ for each $x\in(2\ell/3,\ell)$.
Such choice can be considered as reasonable concerning the hypothesis~\eqref{decay}, and the results of the computation are shown in Figure~\ref{fig:Random1tau}.

\begin{figure}[hbt]
\begin{center}
{\includegraphics[width=6.75cm]{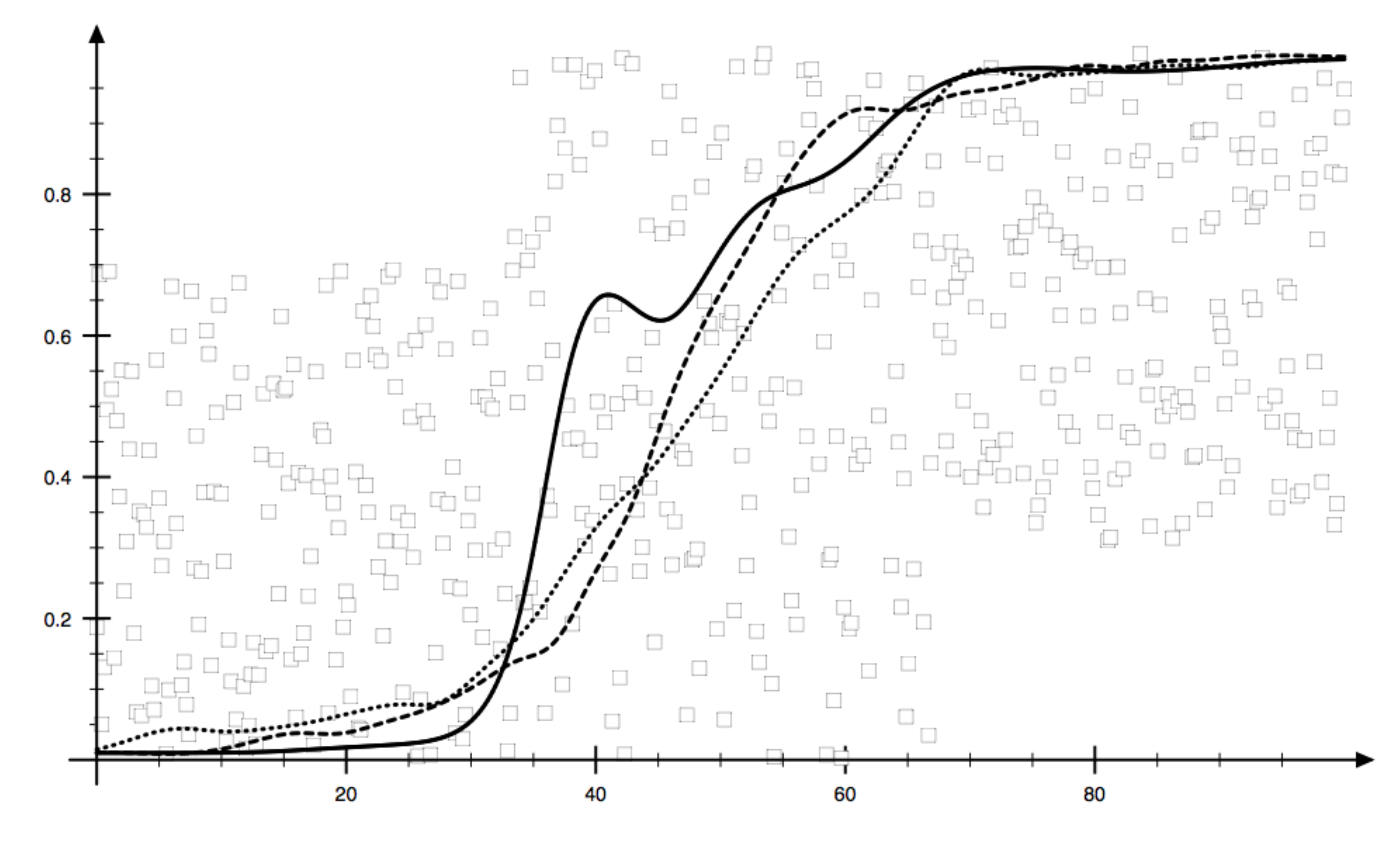}}
{\includegraphics[width=6.75cm]{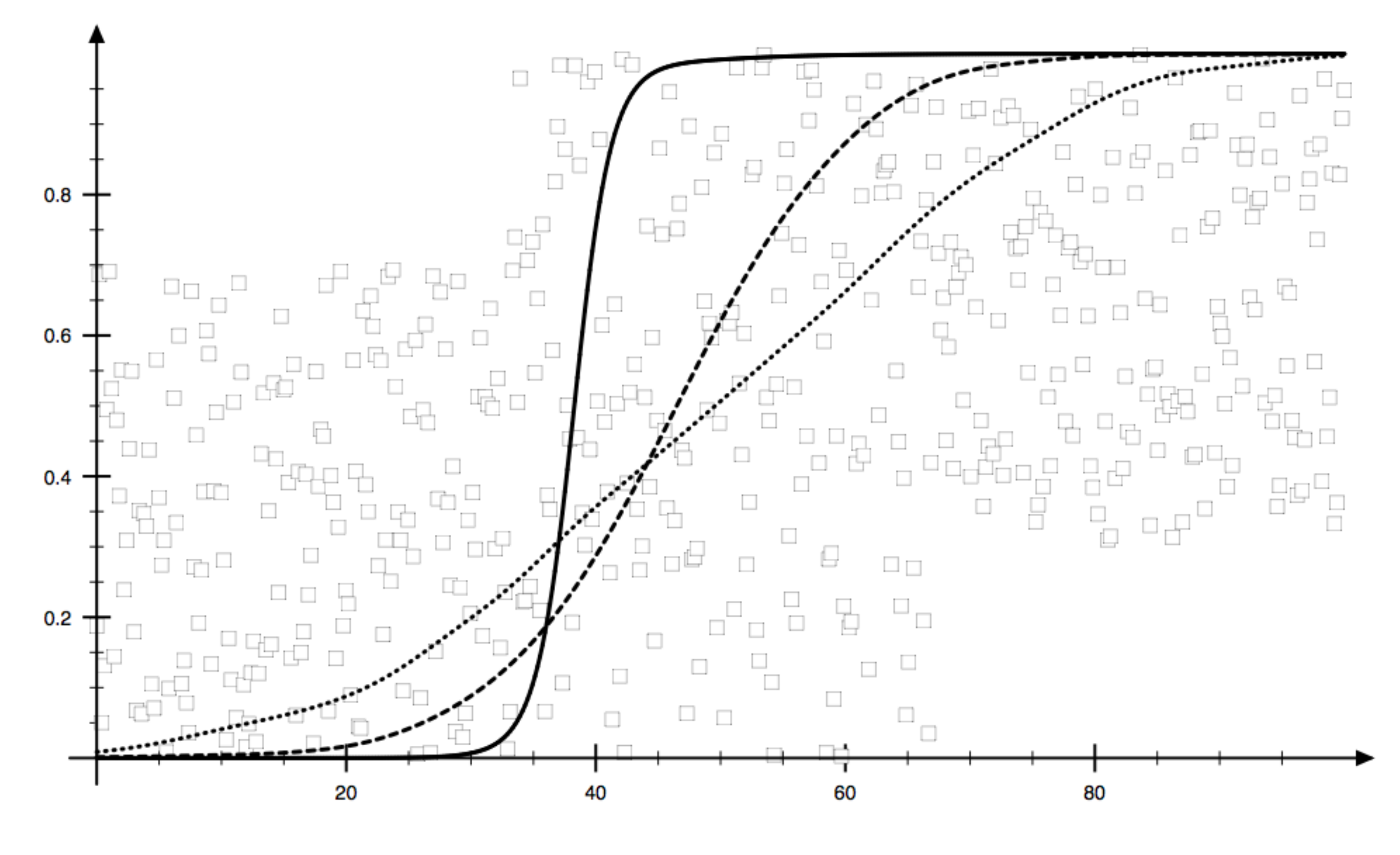}}
\end{center}
\caption{Random initial datum ($\square$) in $(0,\ell)$.
Solution profiles for the hyperbolic Allen--Cahn equation with relaxation
at times $t=10$ (left) and $t=20$ (right),
for $\tau=1$ (continuous line), $\tau=5$ (dashed) and $\tau=10$ (dotted).
}\label{fig:Random2tau}
\end{figure}

The transition is even more robust than what the previous computation shows, since initial data
that do not satisfy the requirement~\eqref{decay} still exhibits convergence.
As an example, let us consider the case of a randomly chosen initial datum $u_0(x)$ given by
any random value in $(0,0.7)$ for each $x\in(0,\ell/3)$, in $(0,1)$ for each $x\in(\ell/3,2\ell/3)$ and in $(0.3,1)$ for each $x\in(2\ell/3,\ell)$.
Also in such a case, we clearly observe the appearance and formation of a stable front, as shown in Figure~\ref{fig:Random2tau}.

\begin{figure}[hbt]
\begin{center}
{\includegraphics[width=6.75cm]{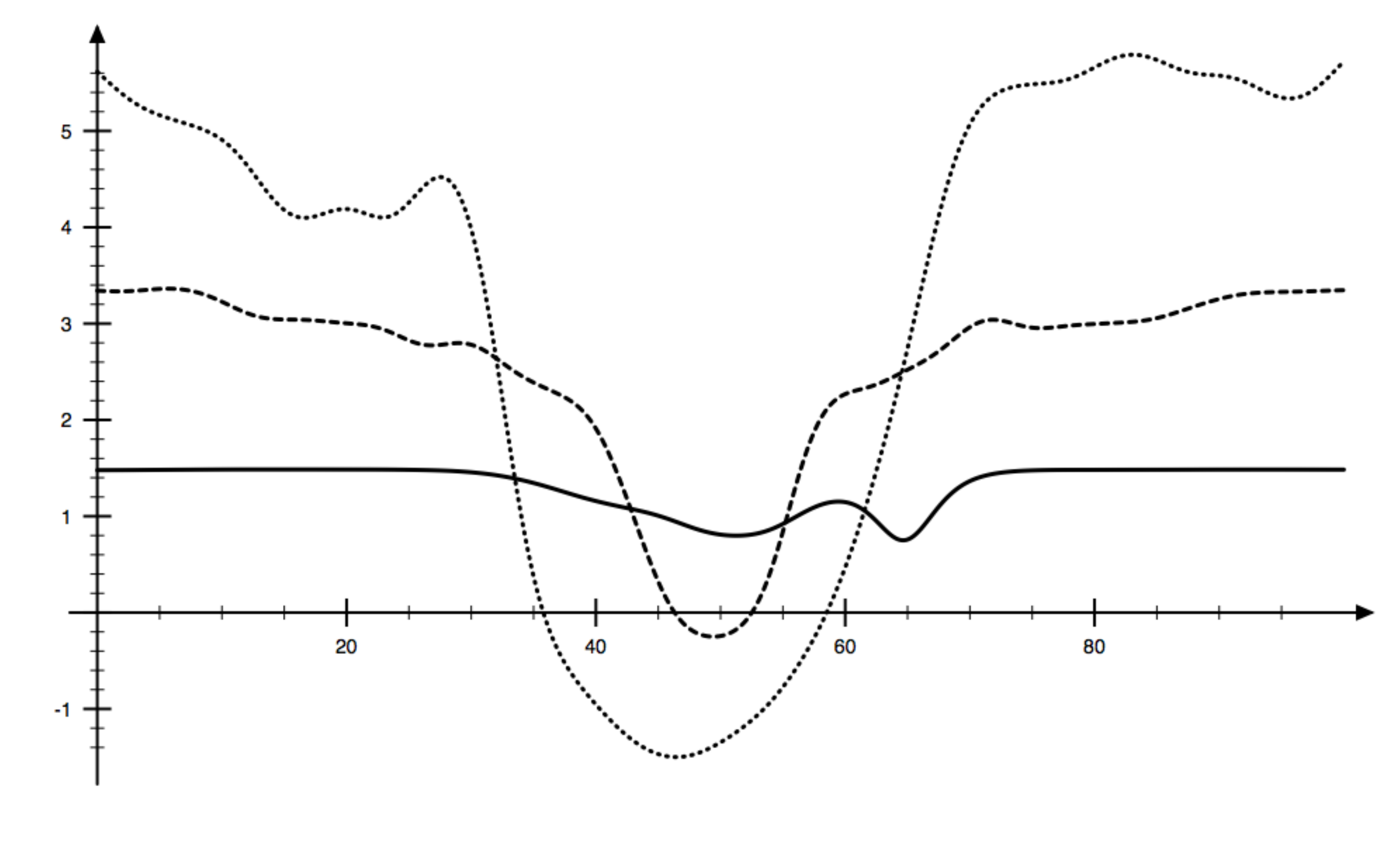}}
{\includegraphics[width=6.75cm]{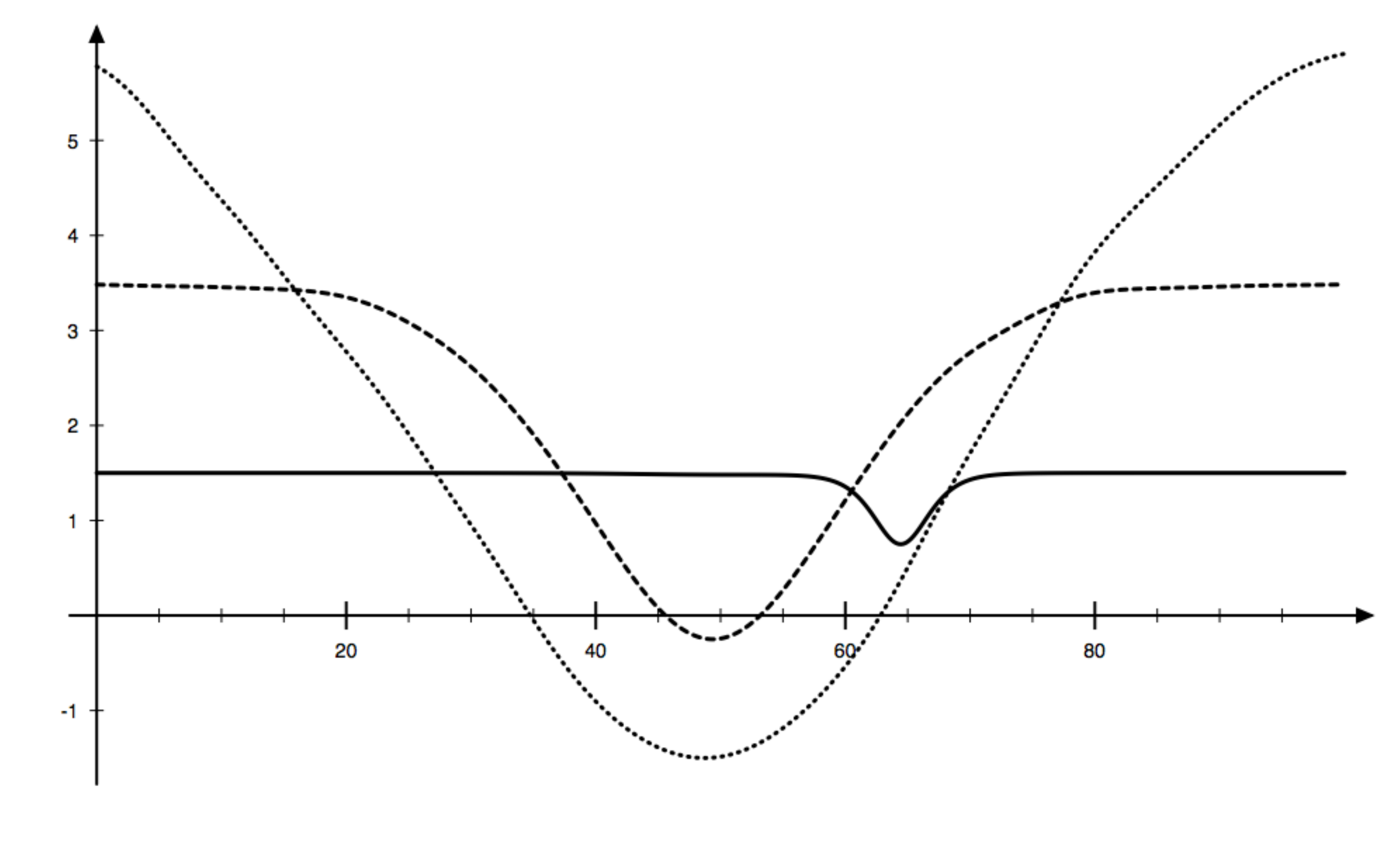}}
\end{center}
\caption{Profile of the function $g(u):=1-\tau f'(u)$ for time $t=10$ (left) and $t=20$ (right)
corresponding to the initial datum shown in Figure~\ref{fig:Random1tau}.
The legend for the lines is the same as in the previous figures.}
\label{fig:Random1g}
\end{figure}

\begin{figure}[hbt]
\begin{center}
{\includegraphics[width=6.75cm]{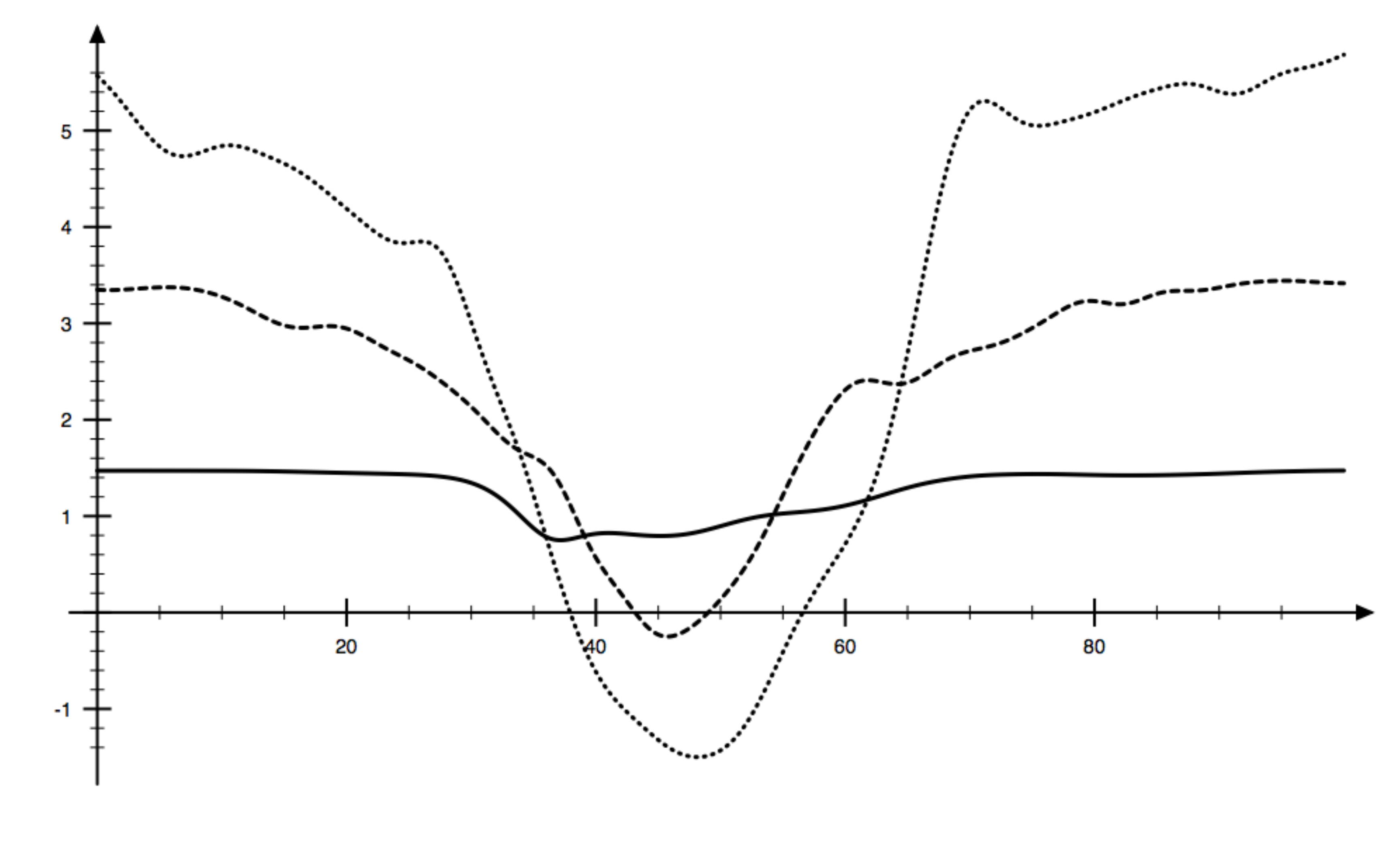}}
{\includegraphics[width=6.75cm]{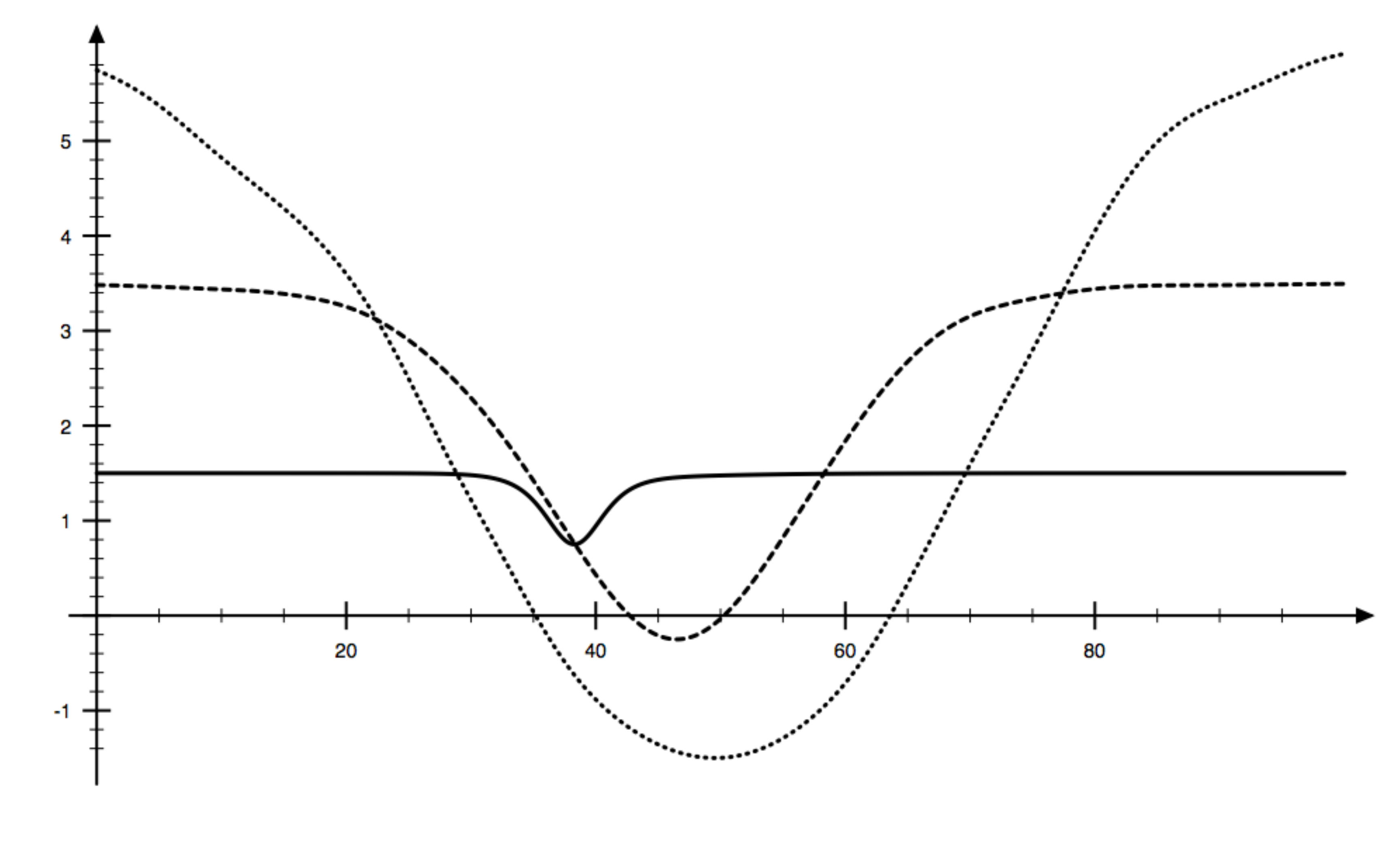}}
\end{center}
\caption{Profile of the function $g(u):=1-\tau f'(u)$ for time $t=10$ (left) and $t=20$ (right)
corresponding to the initial datum shown in Figure~\ref{fig:Random2tau}.
The legend for the lines is the same as in the previous figures.
}\label{fig:Random2g}
\end{figure}

The convergence is manifest also in the case where the stability condition
$g(u):=1-\tau f'(u)>0$ fails in some region.
At least for the cubic (bistable) nonlinear reaction term $f$, such region is typically centered at $u=\alpha$.
In particular, being $u=\alpha$ an unstable equilibrium, $f'(\alpha)$ is positive, thus $g(\alpha)$
is negative when $\tau$ is sufficiently large.
The values of the function $g$ are plotted in Figure~\ref{fig:Random1g} and Figure~\ref{fig:Random2g}, respectively,
for two different times, namely $t=10$ and $t=20$, and different values of $\tau$, namely $\tau=1$, $\tau=5$ and $\tau=10$. Of course, the function $g$ is asymptotically positive, since $0$ and $1$ are stable equilibria, and thus the value of the first order derivative $f'$ is negative. 
The numerical results show that, for sufficiently large values of $\tau\,$, some region corresponding to the center of the wave profile appears where $\tau>1/f'(u)$ for some $u\in(0,1)\,$, and it contains the value $u=\alpha$ (at least for the cubic case).

\subsection{Pseudo-kinetic scheme for the Guyer-Krumhansl's law}

The diagonalization procedure that has been performed in Section~\ref{sec:physics} to deduce a kinetic interpretation of the reaction-diffusion equation with relaxation, starting from the \textit{Maxwell-Cattaneo law}~\eqref{maxwellcattaneo}, it cannot be straightforwardly extended to the case of the \textit{Guyer–Krumhansl law} because of the presence of a higher order (conservative) operator in the model~\eqref{guyerkrumhansl}. Although such an issue is rigorously pursued in a work in progress, here we attempt at presenting an hybrid version of the kinetic scheme~\eqref{semimain} to adapt to the present case, thus providing an easy-to-implement algorithm for the pseudo-parabolic equation~\eqref{guyerkrumRD}.

Starting from~\eqref{semimain}, we consider the following variation,
\begin{equation}
\label{semimain2}
	\frac{d \vv_i}{dt} = - \varrho^2 \frac{\uu_{i+1} - \uu_{i-1}}{2 \dxi} - \frac1{\tau}\,\vv_i 
		+ \bigl( \nu + \frac12 \varrho\,\dxi \bigr) \frac{\vv_{i+1}-2\vv_i+\vv_{i-1}}{\dxi^2}\,,
\end{equation}
that enjoys the same consistency properties as the original scheme, since the order of magnitude of the physical parameter $\nu$ is clearly bigger than that of the correction by the numerical viscosity $\frac12 \varrho\,{\rm dx}\,$.

Although it deserves to be rigorously justified and further confirmed by extensive numerical simulations, this approach is clearly more convenient than the usual way of putting higher order hyperbolic equations like~\eqref{onefield} and~\eqref{guyerkrumRD} in form of lower order systems for numerical issues, namely
\begin{equation*}
	\prt u = w\,, \qquad \tau \prt w + \bigl(1-\tau f^{\prime}(u)\bigr) w - \mu\,\prxx u + \nu\,\prxx w = f(u) + \nu\,\prxx f(u)\,,
\end{equation*}
for which a direct semi-discrete approximation provides $\frac{d \uu_i}{dt} = {\rm w}_i\,$, together with
\begin{equation}
\label{prova1}
	\begin{split}
	\tau \frac{d {\rm w}_i}{dt} = \,& f(\uu_i) - \bigl(1-\tau f^{\prime}(\uu_i)\bigr) {\rm w}_i + \mu\,\frac{\uu_{i+1}-2\,\uu_i+\uu_{i-1}}{{\rm dx}^2}\,, \\
	& - \nu\,\frac{{\rm w}_{i+1}-2\,{\rm w}_i+{\rm w}_{i-1}}{{\rm dx}^2} + \nu\,\frac{{\rm f(u)}_{i+1}-2\,{\rm f(u)}_i+{\rm f(u)}_{i-1}}{{\rm dx}^2}\,.
	\end{split}
\end{equation}

Another way of dealing with higher order one-field equations can be the following: we rewrite~\eqref{guyerkrumRD} as
\begin{equation*}
	\prt \bigl( \tau \prt u + u - \tau f(u) - \nu\,\prxx u \bigr) - \mu\,\prxx u = f(u) - \nu\,\prxx f(u)\,,
\end{equation*}
for which an alternative representation as second order system is given by
\begin{equation*}
	\tau \prt u + u - \tau f(u) - \nu\,\prxx u = w\,, \qquad \prt w - \mu\,\prxx u = f(u) - \nu\,\prxx f(u)\,,
\end{equation*}
thus generalizing~\eqref{reactiondiffusion}, with corresponding semi-discrete approximation
\begin{equation}
\label{prova2}
	\begin{split}
	\tau \frac{d \uu_i}{dt} &= {\rm w}_i - \uu_i + \tau f(\uu_i) + \nu\,\frac{\uu_{i+1}-2\,\uu_i+\uu_{i-1}}{{\rm dx}^2}\,, \\
	\frac{d {\rm w}_i}{dt} &= f(\uu_i) + \mu\,\frac{\uu_{i+1}-2\,\uu_i+\uu_{i-1}}{{\rm dx}^2} - \nu\,\frac{{\rm f(u)}_{i+1}-2\,{\rm f(u)}_i+{\rm f(u)}_{i-1}}{{\rm dx}^2}\,.
	\end{split}
\end{equation}

Both schemes~\eqref{prova1} and~\eqref{prova2} formally converge to the standard discretization of~\eqref{reactiondiffusion} for $\tau\to 0^+$ and $\nu \to 0\,$, but they exhibit the well-known criticality of defining the correct reconstruction of the external field $f(u)$ on the (possibly nonuniform) spatial mesh. Therefore, the pseudo-kinetic scheme~\eqref{semimain2} maintains a wider interest in view of its underlying physical interpretation.

We conclude by remarking that such peculiar feature is not shared by other more general forms of relaxation system, for instance
\begin{equation}
\label{onefield2}
	\tau \prtt u + g(t,x,u\,;\tau)\,\prt u - \mu\,\prxx u = f(u)\,,
\end{equation}
that is considered in~\cite{FLM}, for example. Unless specific expression for the external field $g$ are taken into account for physical reasons, the only approach to the numerical approximation of~\eqref{onefield2} seems to be the transcription into a first order system by putting
\begin{equation*}
	\prt u = w\,, \qquad \tau \prt w + g(t,x,u\,;\tau)\,w - \mu\,\prxx u = f(u)\,.
\end{equation*}
On the other hand, under the hypothesis that $g$ does not depend explicitly on the independent variables, one can consider
\begin{equation*}
	g(u\,;\tau)\,\prt u = \prt \bigl(g(u\,;\tau) u\bigr) - \partial_{u} g(u\,;\tau)\,u\,\prt u
\end{equation*}
and then equation~\eqref{onefield2} reads
\begin{equation*}
	\prt \bigl( \tau\,\prt u + g(u\,;\tau) u \bigr) - \partial_{u} g(u\,;\tau)\,u\,\prt u - \mu\,\prxx u = f(u)\,,
\end{equation*}
so that we can define
\begin{equation*}
	\tau \prt u + g(u\,;\tau) u = w\,, \qquad \prt w - \partial_{u} g(u\,;\tau)\,u\,\prt u - \mu\,\prxx u = f(u)\,,
\end{equation*}
with the second equation rewritten like
\begin{equation*}
	\prt w - \frac1{\tau}\partial_{u} g(u\,;\tau)\,u \left(w - g(u\,;\tau)\,u \right) - \mu\,\prxx u = f(u)
\end{equation*}
that is even different from all the previous versions, thus revealing the great advantage of a physical justification for the models at hands, as already suggested in Section~\ref{sec:physics}.

\section*{Acknowledgements}
This work has been partially supported by CONACyT (Mexico) and MIUR (Italy), through the MAE Program for Bilateral Research, grant no. 146529. The work of RGP was partially supported by DGAPA-UNAM, grant IN100318.




\end{document}